\documentclass{amsart}
\usepackage{amssymb}
\usepackage{amscd}
\usepackage{amsmath}
\usepackage{amsfonts}
\usepackage{a4}
\usepackage{amsthm}

\usepackage{pslatex}

\input xy
\xyoption{all}

\newtheorem{thm}{Theorem}[section]
\newtheorem{defi}[thm]{Definition}
\newtheorem{prop}[thm]{Proposition}
\newtheorem{lem}[thm]{Lemma}

\newenvironment{rem}{\noindent \bf Remark : \rm  \noindent}{}
\newenvironment{rems}{\noindent \bf Remarks : \rm  \noindent}{}

\newcommand{\ra}{\rightarrow}
\newcommand{\thra}{\twoheadrightarrow}

\begin{document}

\author{Niels Borne}
\address{
Universit\`a di Bologna \\
Dipartimento di Matematica \\
Piazza di Porta San Donato, 5 \\
I-40126 Bologna \\
Italy}
\email{borne@dm.unibo.it}
\title[A relative Shafarevich theorem]
{A relative Shafarevich theorem}
\maketitle
\tableofcontents

\begin{abstract}
Suppose given a Galois \'etale cover $Y\ra X$ of proper non-singular
curves over
an algebraically closed field $k$ of prime characteristic $p$.
Let $H$ be its Galois group, and $G$
any finite extension of $H$ by a $p$-group $P$. We give necessary and
sufficient conditions on $G$ to be the Galois group of an \'etale
cover of $X$ dominating $Y\ra X$. 
\end{abstract}

\section{Introduction}
\subsection{Main statement}

Let $p$ be a prime integer, and $k$ an algebraically closed field of
characteristic $p$.
In \cite{PS}, A.Pacheco and K.Stevenson have given necessary and
sufficient conditions for an extension $G$ of a
$p'$-group by a $p$-group to occur as the Galois group
of an \'etale cover of
a given proper, non-singular curve $X$ over $k$.
We pursue this work further by solving the case of the
extension $G$ of an arbitrary group $H$ (which we suppose of course
realized as a Galois group) by a $p$-group $P$.

Thanks to the work of
S.Nakajima (see \cite{Naka}), we can define, for any \'etale
$H$-cover $Y \thra
X$, a family of integers $\delta_{Y,V}$ parametrized by the set
$\mathcal S_H$ of
isomorphism classes of simple $k[H]$-modules. These integers describe
the Galois module structure of the space $H^0(Y, \Omega_Y)^s$
of semi-simple differentials
on $Y$ in the sense that they are defined by the formula :

$$H^0(Y, \Omega_Y)^s \simeq  \Omega^2 k \oplus \mathop{\oplus}_{V \in \mathcal
  S_H} P(V) ^{\oplus \delta_{Y,V}}$$

\noindent (see \S \ref{def} for a precise definition).

These ``Hasse-Witt coefficients'' seem to be slightly more convenient
invariants for our purposes than the previously used
generalized Hasse-Witt invariants (see Remark following definition \ref{HW}).
Let $Z \ra X$ be an \'etale $G$-cover dominating the \'etale $H$-cover   
$Y \thra X$, and $q:G \thra H$ the corresponding
group extension. We show in part \ref{HWcoeff} the following formula
relating the corresponding Hasse-Witt coefficients :

$$ \forall V\in \mathcal S_H \;\;\;\;\; \delta_{Z,q^*V} + \dim_k H^1(G,q^*V) =
\delta_{Y,V} + \dim_k H^1(H,V)$$

\noindent (see proposition \ref{quot}). 

We then use this result in part \ref{extcov} to prove an
extension of the main Theorem of \cite{PS} :

\begin{thm}

  \label{main}

  Let $k$ be an algebraically closed field $k$ of positive characteristic $p$.
Let also $Y \ra X$ be a Galois \'etale cover of proper non-singular curves
over $k$, with Galois group $H$.
Suppose
$$\xymatrix{ 1 \ar[r] & P \ar[r] & G \ar[r]^{q} &  H \ar[r]  & 1  } $$
is an exact sequence of finite groups, where $P$ is a $p$-group.

Then the embedding problem

$$\xymatrix{ &&& \pi_1(X) \ar@{->>}[d]\ar@{-->}[dl] \\1 \ar[r] & P \ar[r] & G \ar[r]^{q} &  H \ar[r]  & 1  }  $$     

\noindent has a strong solution if and only if for every simple $k[H]$-module
$V$ the following inequality holds :

$$\dim_k H^1(G,q^*V)-\dim_k H^1(H,V) \leq \delta_{Y,V}$$

\end{thm}

\ \\
Here are some indications about the origins of this theorem.

When $H=1$, the above result specialises to the classical Theorem of
Shafarevich (see for instance \cite{Crew}) :
a given $p$-group $P$ occurs as the Galois group of an \'etale
cover of $Y$ if and only if the minimal number of generators of $P$ is
less than the Hasse-Witt invariant of $Y$.

The first to consider inequalities as in the Theorem seems to be E.Kani in
\cite{Kani}, and for this reason we call them {\it the Kani
inequalities} or {\it the Kani conditions}.
These were given as necessary conditions for the
existence of a tame cover, and since Kani considered differentials instead
of semi-simple ones, his inequalities coincide with ours
only when the curve $Y$ is ordinary.  

Our main source of inspiration was \cite{PS}, which we tried to
interpret from the ``Galois module'' point of view. The original
problem and the presentation in terms of an embedding problem is
due to the authors of that paper. However, the tools
employed here are quite different, in particular since we use some
basic results of
modular representation theory. Moreover, we make a heavy use of the
convenient Proposition \ref{quot}, which seems new, and of independant
interest. As an application, we define, in paragraph \ref{Vord},
a refinement of the
notion of ordinariness of a Galois \'etale cover, which we call
$V$-ordinariness.

In a last paragraph, we give another proof of Theorem \ref{main}
in the special case solved in \cite{PS}, namely when $H$ is
supposed a $p'$-group (see Proposition \ref{alg}). This proof is based
on the same cohomological argument as the one given in {\cite{PS}},
but it uses moreover the existence of an interesting quotient of the
fundamental group of~$X$.

\subsection{Acknowledgements}

I would like to thank both M.Emsalem and B.Erez for their
support.

\newpage

\section{Hasse-Witt coefficients}

\label{HWcoeff}
\subsection{Notations}

\subsubsection{Curves}

In the following we fix an algebraically closed field $k$ of positive
characteristic $p$.

The letters $X$, $Y$, ... will denote connected, complete,
non-singular curves over $k$.

Let $S$ be an effective divisor on a curve $X$. Then the Cartier operator
$\mathcal C$ acts on the space of global sections of the logarithmic
differentials $H^0(X, \Omega_X (S))$ (one can define this action
directly or as
dual of the action of the Frobenius, see \cite{Ser}). This action is
$p^{-1}$-linear. We will call semi-simple logarithmic differentials
the subspace $$H^0(X, \Omega_X (S))^s :=k
\mathop{\otimes} _{\mathbb F_p} H^0(X, \Omega_X (S))^{\mathcal C}$$
generated by forms fixed by $\mathcal C$ of the whole space $H^0(X,
\Omega_X (S))$ of logarithmic differentials. When $S=0$, this
subspace has dimension $\gamma_X$, the classical Hasse-Witt invariant
of the curve $X$. One defines also a natural supplement of $H^0(X,
\Omega_X (S))^s$ in $H^0(X,\Omega_X (S))$, the space $H^0(X,\Omega_X (S))^n$ of
nilpotent differentials (see \cite{Ser}).

We will denote by $\pi_1(X)$ ``the'' fundamental group of the curve
$X$, without specifying a base point. When considering morphisms
between curves and their fundamental groups, we will always suppose
that a coherent choice of base points is made. 

\subsubsection{Representation theory}

We also need a few basic notions from module theory over an artinian
(non-commutative) ring $A$. All $A$-modules will be implicitly
supposed of finite type. We will use these notions for the group
ring $A=k[H]$ (sometimes $A=\mathbb F_p[H]$) where $H$ is a finite
group whose order is (possibly) divisible by $p$. We refer to
\cite{Bens} for more details.

The letter $k$ will always denote the trivial representation.

Every finitely generated $A$-module $M$ has a well-defined (i.e. unique
up to isomorphism) projective cover, which we denote by $P(M)
\twoheadrightarrow M$. This allows us to define $\Omega M := \ker( P(M)
\twoheadrightarrow M)$ and inductively $\Omega^{i+1} M=\Omega(\Omega^{i} M)$.

When there are several groups acting, we will also use the notations
with indices $P_H(M)$ and $\Omega_H M$ to stress the fact that $H$ acts.

We will denote by $\mathcal S _A$ (or $\mathcal S_H$ when there is no
doubt about the field used) a fixed complete set of
representatives of the isomorphism classes of simple $A$-modules.

\subsubsection{Warning}
We make use of the usual notations and vocabulary of both domains, but this is
not without certain ambiguities. For an example, the letter $\Omega$
is used with two different meanings. We must also point out that the
word semi-simple is very ambiguous. For an example we will see that, if
$Y$ is a curve acted on by a group $H$, the space of semi-simple
differentials on $Y$ is in general not semi-simple as a $k[H]$-module.
We hope that this warning and the context will help to give a clear
meaning to our assertions.

\subsubsection{Profinite groups}
\label{embpb}

An embedding problem is a diagram in the category of profinite groups :
$$\xymatrix{ &&& \pi \ar@{->>}[d]^{a} \\1 \ar[r] & P \ar[r] & G \ar[r]^{q} &  H \ar[r]  & 1  }  $$  
where the vertical arrow is an epimorphism and the horizontal sequence
is exact. It is said to have a {\it weak solution} if there exists a
continuous homomorphism
$\beta: \pi\ra G$ lifting $\alpha$, i.e. $q \circ \beta
=\alpha$. There is a {\it strong solution} if one can choose moreover
$\beta$ to be an epimorphism.
According to \cite{Ser2}, Chapitre I, \S 3.4, Proposition 16, if ${\rm
cd}_p \pi \leq 1$ and $P$ is a pro-$p$-group, there is always a weak
solution.

\subsection{Definition}
\label{def}

\begin{lem}
Let $Y \ra X$ be a Galois \'etale cover of proper non-singular curves,
with Galois group $H$. Then

(i) For any reduced non-trivial effective divisor $S$ the space
$H^0(Y, \Omega_Y
(S))^s$ of semi-simple differentials logarithmic along $S$ is
projective as a $k[H]$-module.

(ii) There is a projective $k[H]$-module $P$, unique up to isomorphism,
such that $H^0(Y, \Omega_Y)^s \simeq  \Omega^2 k \oplus P$
\end{lem}

\begin{proof}

(i) Let $P$ be a $p$-Sylow of $H$. According to \cite{Naka}, Theorem 1,
applied to the subcover $Y \ra Y/P$, the space $H^0(Y, \Omega_Y
(S))^s$ is $k[P]$-free. But this is equivalent to say that this space
is $k[H]$-projective (this is easily seen using the fact that any
$k[H]$-module is projective relative to $P$ [see for an example
\cite{Bens} for relative projectivity]).

(ii) Choose $S$ to be reduced to an orbit under $H$. Then the short exact
sequence of sheaves $ 0 \ra \Omega_Y \ra \Omega_Y(S) \ra \mathcal
O_S\ra 0$
gives the long exact sequence

$$ 0 \ra H^0(Y, \Omega_Y)^s \ra H^0(Y, \Omega_Y(S))^s \ra k[H] \ra k \ra
0$$

\noindent (see the proof of \cite{Naka}, Theorem 2 for the fact that nilpotent
differentials have trivial residues).
In view of (i), it suffices to apply the well-known algebraic lemma
comparing this sequence to the beginning of the minimal projective
resolution of $k$, see for instance \cite{Kani}, Proposition 4.

\end{proof}

Now, using the Krull-Schmidt theorem to decompose $P$, and the fact
that the indecomposable projective $k[H]$-modules are exactly the
projective covers of the simple $k[H]$-modules, we can make the structure
of $H^0(Y, \Omega_Y)^s$ more precise :

\begin{defi}
\label{HW}
Let $Y \ra X$ be a Galois \'etale cover with Galois group $H$ and let
$\mathcal S_H$ be a set of representatives of the isomorphism classes
of simple $k[H]$-modules. We define the
{\em Hasse-Witt coefficients} of the cover as the only non-negative integers
$\delta_{Y,V}$ verifying~:
$$H^0(Y, \Omega_Y)^s \simeq  \Omega^2 k \oplus \mathop{\oplus}_{V \in \mathcal
  S_H} P(V) ^{\oplus \delta_{Y,V}}$$

\end{defi}

\begin{rem}
It would be more coherent with previous work on the subject to
consider generalized Hasse-Witt invariants defined by 
$\gamma_{Y,V}:=\delta_{Y,V}\cdot \dim_k P(V)$. However, since $\dim_k
P(V)$ is hard to compute in general, and contains no geometric
information, we prefer to exclude it from the definition,
to concentrate on the computation of $\delta_{Y,V}$'s. 
\end{rem}

\subsection{Behaviour under quotient : statement}

First we introduce a useful notation.

\begin{defi}
If $M$ is any $k[H]$-module of finite type, we define
$$h^i(H,M):=\dim_k H^i(H,M)$$
\end{defi}

If $q: G \thra H$ is a group epimorphism and $M$ is a $k[H]$-module, we
denote by $q^{*}M$ the corresponding inflated $k[G]$-module. (We will also
use $M=q_* N$ for $N=q^{*}M$).

\begin{prop}
\label{quot}
Let $Z \ra X$ be a Galois \'etale cover with Galois group $G$, and let
$Y \ra X$ be any quotient Galois cover, with Galois group $H$.
Denote by $q$ the quotient map $G \twoheadrightarrow H$. Then :

$$ \forall V\in \mathcal S_H \;\;\;\;\; \delta_{Z,q^*V} + h^1(G,q^*V) =
\delta_{Y,V} + h^1(H,V)$$

\end{prop}

Before giving the proof of Proposition \ref{quot}, we make some
remarks, beginning by another useful definition.

\begin{defi}
\label{KaniCond}
Let $Y \ra X$ be a Galois $H$-cover, $\pi _1(X) \twoheadrightarrow H$ the
corresponding group epimorphism, and $q : G \twoheadrightarrow H$
an extension of $H$. We will say that the pair $( \pi
_1(X) \twoheadrightarrow H, G \twoheadrightarrow H)$  satisfies
  the Kani inequalities (or conditions) if and only if :
$$ \forall V\in \mathcal S_H \;\;\;\;\;  h^1(G,q^*V)- h^1(H,V) \leq \delta_{Y,V}$$
\end{defi}

\newpage
\begin{rems}
  \begin{enumerate}
\item  In the perspective of the proof of Theorem \ref{main}, it is clear
  that Proposition \ref{quot} implies that if the embedding
  problem of the theorem has a solution, then the Kani
  inequalities hold for $( \pi
_1(X) \twoheadrightarrow H , G \twoheadrightarrow H)$. 
  
\item Proposition \ref{quot} is valid without assuming that the
  kernel $P$ of $q : G \twoheadrightarrow H$ is a $p$-group. However, this
  hypothesis is sometimes useful. For an example, it implies that the
  inflation map $q^*: \mathcal S_H \ra \mathcal S_G$ is one to
one. Indeed, we have :
\begin{thm}[Clifford]
\label{clif}
If $P$ is any normal subgroup of $G$ and
$M$ is a semi-simple $k[G]$-module, then $M_{|P}$ is
a semi-simple $k[P]$-module. 
\end{thm}
\noindent (see for instance \cite{Alp}, Chapter I.3, 
Theorem 4). Moreover, it is well known that if $P$ is a $p$-group, the
only simple $k[P]$-module is the trivial one.
Hence, if we suppose that $P$ is a $p$-group, any
simple $k[G]$-module is trivial when restricted to $P$.

We can conclude, with this hypothesis on $P$,
that the Hasse-Witt coefficients of the covers $Z \ra X$ and $Y \ra
X$ determine each other.

\item M.Emsalem suggested to us that one can consider more generally
tame covers, and that Proposition \ref{quot} could also hold in
this more general setting. We hope to come back to these questions in a
subsequent paper. 

\end{enumerate}

\end{rems}

\subsection{Behaviour under quotient : proof}

The idea is quite naturally to take the fixed part under the kernel $P$ of
$q : G \twoheadrightarrow H$ of the two members of the isomorphism~:
$$H^0(Z, \Omega_Z)^s \simeq  \Omega^2_G k \oplus \mathop{\oplus}_{U
  \in \mathcal S_G} P_G(U) ^{\oplus \delta_{Z,U}}$$

\noindent and then use the three following lemmas.

\subsubsection{Algebraic part}
In this part we assume only that we have an arbitrary exact sequence
of groups :

$$\xymatrix{ 1 \ar[r] & P \ar[r]^{r} & G \ar[r]^{q} &  H \ar[r]  & 1  } $$

\begin{lem}
  \label{lem1}
Let $U\in \mathcal S_G$. $$q_*(P_G(U)^P) \simeq  \left\{ \begin{array}{cl}
P_H(V)& {\rm if}\;\; \exists \;\;V \in \mathcal S_H \;\;{\rm such
  \;that} \;\; U=q^*V \\
  0       & {\rm  else}
    \end{array} \right .$$
\end{lem}
\begin{proof}
The proof is a straightforward computation, using the Krull-Schimdt
theorem and three facts :
\begin{enumerate}
\item If $U \in \mathcal S_G$, then ${\rm soc}(P_G(U)) \simeq U$,
  \newline (see \cite{Alp}, Chapter II.6, Theorem 6, or \cite{Bens}, Theorem 1.6.3)~;
\item $q_*(k[G]^P) \simeq k[H]$, \newline so in particular $q_*(\cdot^P )$
  sends projectives to projectives ;
\item $k[G] \simeq \mathop{\oplus}_{U
  \in \mathcal S_G} P_G(U) ^{\oplus \dim_k U}$.
\end{enumerate}

\end{proof}

\begin{lem}
$$q_*((\Omega^2_G k)^P) \simeq \Omega^2_H k\oplus \mathop{\oplus}_{V
  \in \mathcal S_H} P_H(V)^{\oplus h^1(G,q^*V)- h^1(H,V)}$$
\end{lem}

\begin{proof}
  We will suppose that $p | \#P$, the proof in the opposite case being
  similar, and in fact easier. Then, as Kani notices (\cite{Kani}, \S
  5, Lemma 2), one can easily show inductively that for all $i$,
  the module $\Omega^i_P k$ is non-projective, and this implies $H^1(P,
  \Omega^i_P k) \neq 0$ (see \cite{Bens}, Proposition 3.14.4).

  Restricting the exact sequence
  $$\xymatrix{ 0 \ar[r] & \Omega^i_G k  \ar[r] & P_G(\Omega^{i-1}_G k)  \ar[r]
    & \cdots \ar[r] & P_G(\Omega_G k)  \ar[r]  & k \ar[r] & 0   }
  $$

  \noindent
  to $P$ via $r: P \ra G$ (which, in view of $r^*k[G] \simeq k[P]^{\oplus \#
  H}$, also sends projectives to projectives) and comparing to the
  minimal $k[P]$-resolution of the trivial module $k$, we get an
  isomorphism $r^*(\Omega^i_G k) \simeq \Omega^i_P k \oplus Q$, where
  $Q$ is a projective $k[P]$-module, and this implies
  $H^1(P, r^*(\Omega^i_G k) )\simeq H^1(P, \Omega^i_P k)$. So we have
  shown in particular that for any $i$, we have $H^1(P, r^*(\Omega^i_G k))
  \neq 0$.

  Now applying the functor $q_*(\cdot ^P )$ to the exact sequence
  $$\xymatrix{ 0 \ar[r] & \Omega_G k \ar[r] &  P_G k  \ar[r] &  k \ar[r] & 0}$$ we get
  the exact sequence of $k[H]$-modules :
  $$\xymatrix{ 0 \ar[r] & q_*((\Omega_G k) ^P ) \ar[r] & q_*(( P_G k)
    ^P ) \ar[r] &  k \ar[r] &  H^1(P,\Omega_G k)  \ar[r] & 0}$$

  The $k$-vector-space underlying $H^1(P,\Omega_G k)$ is $H^1(P,
  r^*(\Omega_G k))\neq 0$, hence the middle arrow has to be zero, and
  thanks to Lemma \ref{lem1} it comes :
  $q_*((\Omega_G k) ^P )\simeq q_*(( P_G k)^P ) \simeq P_H(k)$.

  We proceed similarly, starting with the sequence 

  $$\xymatrix{ 0 \ar[r] & \Omega^2_G k \ar[r] &  P_G (\Omega_G k)
    \ar[r] & \Omega_G k \ar[r] & 0}$$
we get the exact sequence of $k[H]$-modules :

$$  \xymatrix{ 0 \ar[r] & q_*((\Omega^2_G k) ^P ) \ar[r] & q_*((
    P_G (\Omega_G k))
    ^P ) \ar[r] &  P_H(k) \ar[r] &  H^1(P,\Omega^2_G k)  \ar[r] & 0}
  $$

  We claim that $H^1(P,\Omega^2_G k) \simeq k$. In fact the underlying
  $k$-vector-space is \newline $H^1(P,r^*\Omega^2_G k)\simeq H^1(P, \Omega^2_P
  k)$. But the above exact sequence for $G=P$ (hence $H=1$) shows that
  $H^1(P, \Omega^2_P k) \simeq k$. So $H^1(P,\Omega^2_G k)$ must be
  simple, and because there is a surjective $k[H]$-map $P_H(k)
  \twoheadrightarrow H^1(P,\Omega^2_G k)$, it can only be the trivial
  simple module.
  From this we deduce, comparing the above exact sequence to the
  minimal $k[H]$-projective resolution of $k$, that $ q_*((\Omega^2_G k)
  ^P) \simeq \Omega^2_H k \oplus R$, where $R$ is a projective
  $k[H]$-module. Now, using Lemma \ref{lem1}, and the isomorphism

  $$P_G(\Omega_G k) \simeq \mathop{\oplus}_{U \in \mathcal S_G} P_H(U)^{\oplus
    h^1(G,U)}$$

    \noindent (see \cite{Kani}, Proposition 5), a direct computation in
    the Grothendieck ring of the finitely-generated $k[H]$-modules
    shows the egality between classes :

    $$ [R]= [\mathop{\oplus}_{V \in \mathcal S_H} P_H(V)^{\oplus
    h^1(G,q^*V)- h^1(H,V)}]$$

  \noindent But since both modules involved are $k[H]$-projective, they have to
  be isomorphic.

\end{proof}

\subsubsection{Geometric part}

In this part the hypothesis are those of Proposition \ref{quot}.

\begin{lem}
The pullback along $ Z \ra Y$ induces isomorphisms

(i) of $\mathbb F _p [H]$-modules :
$$ {\rm Pic}_0(Y)[p] \simeq ({\rm Pic}_0(Z)[p])^P $$

(ii) of $k[H]$-modules :
$$ H^0(Y, \Omega_Y)^s \simeq  (H^0(Z, \Omega_Z)^s)^P $$
\end{lem}

\begin{proof}

(i)
Since the map $ Z \ra Y$ commutes with the action of $G$, it is clear
that the pullback is equivariant.

The Hochschild-Serre spectral sequence (see for instance
\cite{Milne}, Chapter III, Theorem 2.20) for the $P$-Galois cover $ Z
\ra Y$ and the sheaf $ \mathbb G_{m,Y}=\mathcal O _Y ^*$ gives the
exact sequence of terms of low degree :

$$\xymatrix{ 0 \ar[r] & H^1(P,k^*) \ar[r] & {\rm Pic}(Y)
     \ar[r] & ({\rm Pic}(Z))^P \ar[r] &  H^2(P,k^*)  }$$

Now since $k$ is of characteristic $p$, raising to the power $p$ induces
a group isomorphism

\xymatrix@R=1pt{ &&&&& k^* \ar[r] & k^*  \\&&&&& x \ar@{|->}[r] & x^p
  }

In particular, the groups $H^i(P,k^*)$ have no $p$-torsion and are
$p$-divisible , and an immediate diagram chase allows to conclude.

\ \\
(ii) According to \cite{Ser}, Proposition 10, there is a natural
isomorphism
$$ {\rm Pic}_0(Y) [p] \simeq H^0(Y, \Omega_Y)^{\mathcal C}$$ 
Indeed, if $\mathcal L$ is an invertible sheaf on $Y$ such that $ \mathcal L
^{\otimes p} \simeq \mathcal O_X$, write $\mathcal L \simeq \mathcal
L_Y(D)$ for a divisor $D$ on $Y$, and let $f$ be a function on $Y$
such that $pD = (f)$. To $\mathcal L$ we can then associate the form
$\omega = df/f$, and this gives the map we were looking for.
Thanks to this description, one checks immediatly that this map is
compatible with the action of $H$, and also with the pull-back along $Z
\ra Y$. But then (i) gives

$$  H^0(Y, \Omega_Y)^{\mathcal C} \simeq (H^0(Z, \Omega_Z)^{\mathcal
  C})^P$$

\noindent and the lemma follows.

\end{proof}

\subsection{$V$-ordinariness}
\label{Vord}
This section is not necessary for the proof of Theorem \ref{main}
but is just an illustration of the usefulness of Proposition
\ref{quot}. Recall that a cover $Y \ra X$ is said ordinary if the
curve $Y$ is, that is if $\gamma_Y=g_Y$. 

\begin{defi}
Let $Y \ra X$ be an \'etale Galois cover of group $H$ and $V \in
\mathcal S_H$. 

The cover is said $V$-ordinary if $\delta_{Y,V}=\dim V(g_X -1)-h^1(H,V)+h^0(H,V)$.
\end{defi}

The link with the classical notion is given by the following lemma.

\begin{lem}
The cover $Y \ra X$ is ordinary if and only if it is $V$-ordinary for
all $V \in \mathcal S_H$.
\end{lem}

\begin{proof}
$Y \ra X$ is ordinary if and only if
$H^0(Y,\Omega_Y)^s=H^0(Y,\Omega_Y)$.
But \cite{Kani}, Theorem 2, gives that
$$H^0(Y,\Omega_Y) \simeq 
\Omega^2 k \oplus \mathop{\oplus}_{V \in \mathcal
  S_H} P(V) ^{\oplus \dim V(g_X -1)-h^1(H,V)+h^0(H,V)}$$
and the lemma follows.
\end{proof}

It is well known that a Galois \'etale $p$-cover of an ordinary curve is
ordinary (this results from the comparison of the Hurwitz and the
Shafarevich formulas). Here is a relative version :

\begin{lem}
Let $Z\ra X$ and $Y \ra X$ be covers as in Proposition \ref{quot}.
Then $Z\ra X$ is $q^{*}V$-ordinary if and only if $Y\ra X$ is $V$-ordinary.
\end{lem}

\begin{proof}
This follows immediately from Proposition \ref{quot}.
\end{proof}

\begin{rem}
This two lemmas suggest a method to show that a cover is
ordinary, by using subcovers.
\end{rem}

\section{Extension of covers}

\label{extcov}

Let $Y \ra X$ be a Galois \'etale cover of Galois group $H$, and
$\pi _1(X) \twoheadrightarrow H$ the corresponding group epimorphism.
In this section, we show that if the pair
$(\pi _1(X) \twoheadrightarrow H, G \twoheadrightarrow H)$
  verifies the Kani inequalities, and the kernel $P$ of
  $G \twoheadrightarrow H$ is a $p$-group, then the
corresponding embedding problem has a solution (see Theorem \ref{main}).

\subsection{Reductions}

In this section, we give an adaptation of an argument of A.Pacheco
and K.Stevenson (see \cite{PS}).

\subsubsection{Reduction lemma}
\begin{lem}

  \label{red}
Let $Y \ra X$ be a Galois \'etale cover of Galois group $H$, and
$\pi _1(X) \twoheadrightarrow H$ the corresponding group epimorphism.
Suppose
$$\xymatrix{ 1 \ar[r] & P \ar[r] & G \ar[r]^q &  H \ar[r]  & 1  } $$
is an exact sequence of finite groups, $P$ is a $p$-group,
and $O$ is a subgroup of $P$ normal in $G$.
Suppose also that the Kani inequalities hold for
$(\pi _1(X) \twoheadrightarrow H, G \twoheadrightarrow H)$.

\ \\
(i) The Kani inequalities hold for $(\pi _1(X) \twoheadrightarrow H, G/O \twoheadrightarrow H)$.

\ \\
(ii) Suppose $\pi _1(X) \twoheadrightarrow G/O$ is a strong solution
of the embedding problem

$$\xymatrix{ &&& \pi_1(X) \ar@{->>}[d]\ar@{-->}[dl] \\1 \ar[r] & P/O
  \ar[r] & G/O  \ar[r]^{q'} &  H \ar[r]  & 1  }  $$

Then the Kani inequalities hold for $(\pi _1(X)
\twoheadrightarrow G/O, G \twoheadrightarrow G/O)$.

\end{lem}

\begin{proof}

We make constant use of the argument in the second remark following
Definition \ref{KaniCond} : in the situation given, Clifford's theorem
(see Theorem \ref{clif}) ensures that the maps $q : G
\twoheadrightarrow H$, $q' : G/O \twoheadrightarrow H$  and $q'':
G \twoheadrightarrow G/O$ induce canonical bijections between the sets
$\mathcal S_H$, $\mathcal S_G$, and $\mathcal S_{G/O}$.

(i) Let $V \in \mathcal S_H$. Note that
$q''_*((q^*V)^O)=q''_*(q^*V)={q'}^* V$. The inflation-restriction sequence
(see for instance \cite{Bens2} 3.5) for $q'': G \twoheadrightarrow
G/O$ and $q^*V$ ensures that $h^1(G/O, {q'}^*V) \leq h^1(G, q^*V)$ and
(i) follows.

(ii) This is a formal consequence of the remark at the beginning of
the proof and of
Proposition \ref{quot}.

\end{proof}

\subsubsection{First reduction}

Suppose that Theorem \ref{main} holds with the extra hypothesis that $P$
is abelian $p$-elementary. Then we claim the theorem holds in
general. This is easily proved by induction on $N=\log_p \#P$, thanks
to Lemma \ref{red}, applied with $O=\Phi(P)$, the Frattini subgroup of
$P$. Notice that $\Phi(P)$ is stable by any automorphism of $P$, hence
in particular normal in $G$.

\ \\
\begin{rem}
The fact that we allow some $P$-part in $H$ enables us to avoid the use
of the argument that the $p$-cohomological dimension of $\pi_1(X)$ is at most
$1$ (compare with \cite{PS}). However, we cannot avoid the use of this
argument at a later stage.
\end{rem}
\subsubsection{Second reduction}

Suppose next that Theorem \ref{main} holds with the extra
hypothesis that $P$ is abelian $p$-elementary, and simple as $\mathbb
F_p[H]$-module. Then we claim the theorem holds with the restricted
extra hypothesis that $P$ is abelian $p$-elementary.
The proof is the same, by induction on $\dim_{\mathbb
  F _p} P$, thanks to Lemma \ref{red}, applied with $O$ the
kernel of $P \twoheadrightarrow Q$, where $Q$ is any $\mathbb F_p[H]$-simple
quotient of $P$. As $H$ acts by conjugation on $P$, the fact that $O$
is $H$-stable ensures that $O$ is normal in $G$.

\subsection{Fundamental exact sequence}

\subsubsection{Introduction}
  
In the following, the data will be :

\begin{enumerate}

\item A Galois $H$-cover $Y \ra X$, or equivalently, a continuous
  surjective group homomorphism $\pi_1(X) \thra H$,

\item A simple $\mathbb F_p[H]$-module of finite type $P$, i.e. an irreducible
  representation of $H$ over $\mathbb F _p$,

\item A group $G$, extension of $H$ by the $\mathbb F_p[H]$-module $P$ (or
  equivalently a class $\alpha \in H^2(H,P)$), such that the
  Kani conditions hold for $(\pi_1(X) \thra H, G \thra H)$.

\end{enumerate}

Of course 2 and 3 are equivalent to the previous data of an exact
sequence $1 \ra P \ra G \ra H \ra 1$, where  $P$ is abelian
$p$-elementary, simple as $\mathbb F_p[H]$-module, such that the
Kani conditions hold for $(\pi_1(X) \thra H, G \thra
H)$.

We split this data in two because 3 appears only at a later stage, and
also because we will have to distinguish two cases, depending on whether
$\alpha$ is zero or not.


\subsubsection{Inflation-restriction sequence for profinite groups}

In the next general proposition, we use bold letters to distinguish
from our previous notations.

\begin{defi}
Let $1 \ra {\bf F} \ra {\bf G} \ra {\bf H} \ra 1$ be an exact sequence of
profinite groups, $M$ a discrete abelian group endowed with a
continuous action of ${\bf G}$. The {\em transgression} morphism is
the fourth arrow in the exact sequence of terms of low degree in the 
Lyndon-Hochschild-Serre spectral sequence :

$$0 \ra H^1({\bf H}, M^{\bf F}) \ra H^1({\bf G}, M) \ra H^1({\bf F}, M)^{\bf H}
\ra H^2({\bf H}, M^{\bf F})\ra H^2({\bf G}, M) $$

\end{defi}

\begin{prop}
\label{infres}
Suppose moreover that ${\bf F}$ acts trivially on $M$. Let ${\bf F'}$ be
the commutator subgroup of ${\bf F}$, and let $\gamma \in H^2({\bf H},
{\bf F}/{\bf F'})$ be the class of the extension :

$$ 1 \ra {\bf F}/{\bf F'} \ra  {\bf G}/{\bf F'} \ra {\bf H} \ra 1$$

Then the transgression ${\rm trans} : H^1({\bf F}, M)^{\bf H}
\ra H^2({\bf H}, M^{\bf F})$ is given explicitly by :
$$ \forall u \in  H^1({\bf F}, M)^{\bf H}=H^0({\bf H}, {\rm Hom}({\bf
  F}/{\bf F'}, M))\;\;\;\; {\rm trans}(u) = -\gamma \cup u $$

\end{prop}

\begin{proof}
We refer to \cite{Hoch}, Theorem 4 for the case of abstract groups.

\end{proof}

\begin{prop}
\label{fond}
 Let $Y \ra X$ be a Galois \'etale cover with Galois group $H$, and $P$
 a simple $\mathbb F_p[H]$-module.

 (i) The exact sequence of terms of low degree for the spectral
 sequence of the extension $1 \ra \pi_1(Y) \ra \pi_1(X) \ra H \ra 1$
 is the following
 exact sequence of finite-dimensional $\mathbb F_p$-vector spaces :

 $ 0 \ra H^1(H, P) \ra H^1(\pi_1(X),P) \ra {\rm Hom}_H(\pi_1^{ab}(Y),P)
 \ra H^2(H,P) \ra 0 $

 (ii) The transgression map ${\rm trans}: {\rm Hom}_H(\pi_1^{ab}(Y),P) \ra
 H^2(H,P)$ in the sequence above can be described as follows : let $u
 \in  {\rm Hom}_H(\pi_1^{ab}(Y),P)$, $u\neq 0$. Then $u$ corresponds
 canonically to a Galois cover $Z\thra X$ whose Galois group is an
 extension of $H$ by $P$, and ${\rm trans}(u)$ is precisely the
 opposite of the class of this extension in $H^2(H,P)$.

\end{prop}

\begin{proof}

 (i) To show the exactness it suffices, in view of Proposition
\ref{infres}, to show that $H^2(\pi_1(X),P)=0$. But $P$ is killed by
$p$, so this equality follows directly from ${\rm
  cd}_p(\pi_1(X))\leq 1$ (for which we refer to \cite{PS}, Theorem 5.11).
 Moreover, the largest $p$-elementary quotient of $\pi_1^{ab}(X)$ is
 of finite dimension over $\mathbb F_p$ (this dimension is nothing
 else than the Hasse-Witt invariant $\gamma_X$), hence so is ${\rm
   Hom}_H(\pi_1^{ab}(Y),P)$  and the $H^i(H,P)$ are also
 known to be finite-dimensional over $\mathbb F_p$. 

(ii) One point here is because $u$ is $H$-equivariant and not equal
$0$, and $P$ is a simple $\mathbb F_p [H]$-module, it has to be
surjective. Hence $u$ determines a $P$-Galois cover $Z \thra Y$, and
because $u$ is $H$-invariant, $\pi_1(Z)$ is in fact
normal in $\pi_1(X)$, which means that the composite cover $Z \thra X$
is Galois. The last assertion can be deduced from Proposition
\ref{infres} by a direct computation, using the explicit
description of cup product of cochains given in \cite{Hoch}, Chapter
II, \S1.
\end{proof}

\subsubsection{Application to the embedding problem}

We see that Proposition
\ref{fond} implies at once that if the group $G$ is a
non-trivial extension of $H$ by $P$ (that is, its class $\alpha$ in
$H^2(H,P)$ is non zero), then the corresponding embedding problem has
a strong solution (i.e. $G$ is the Galois group of a cover of $X$
dominating the $H$-cover $Y \thra X$). One can of course also
see this quite directly, without using the inflation-restriction
sequence, but ``only'' the inequality ${\rm
  cd}_p(\pi_1(X))\leq 1$.
A rather striking fact is that we do not even use the Kani
inegalities here. 

So the only thing which remains to be proved to end the proof
of Theorem \ref{main} is that we can also realize the trivial
extension. The Proposition \ref{fond} implies that it is equivalent to
the following inequality :

$$\dim_{\mathbb F_p}{\rm Hom}_H(\pi_1^{ab}(Y),P) > \dim_{\mathbb F_p}
H^2(H,P) $$

\subsection{End of the proof}

\subsubsection{Rational representations and base extension}

\begin{lem}

\label{rat}

Let $\mathcal V$ a simple $\mathbb F_p [H]$-module.

(i) $k \otimes_{\mathbb F_p} \mathcal V$ is a semi-simple $k[H]$-module.

(ii) $\dim _{\mathbb F_p} H^i(H,\mathcal V)=\dim _k H^i(H,k \otimes_{\mathbb
  F_p} \mathcal V)$

(iii) The multiplicity of a simple $k[H]$-module
$V$ in $k \otimes_{\mathbb F_p} \mathcal V$ is at most 1.

If moreover $V$ and $V'$ are two simple summands of $k \otimes_{\mathbb
  F_p} \mathcal V$ then :

(iv) $h^i(H,V)=h^i(H,V')$,

(v)  for any Galois \'etale $H$-cover $Y\ra X$ we have $\delta_{Y,V}=\delta_{Y,V'}$.

\begin{proof}

(i) There is an obvious epimorphism :

$$ k[H] \thra \frac{k[H]}{ k \otimes_{\mathbb F_p} {\rm rad}\,{\mathbb
    F_p}[H]}
\simeq k \otimes_{\mathbb F_p} \frac{ {\mathbb F_p}[H]}{{\rm rad}\,
  {\mathbb F_p}[H]}$$  and Wedderburn's theorem (\cite{Bens}, Theorem 1.3.5)
ensures that this last algebra is semi-simple. This shows that
${\rm rad} (k[H]) \subset  k \otimes_{\mathbb F_p} {\rm rad}\, {\mathbb
F_ p}[H]$, and because $\mathcal V$ is killed by ${\rm rad}\, {\mathbb
F_ p}[H]$, the module $k \otimes_{\mathbb F_p} \mathcal V$ is killed by
${\rm rad} (k[H])$, hence is semi-simple as a $k[H]$-module.

(ii) This is true without any assumption on $\mathcal V$ because
$$H^i(H,k \otimes_{\mathbb F_p} \mathcal V) \simeq k \otimes_{\mathbb
F_p} H^i(H,\mathcal V )$$

(iii), (iv) The proofs are rather technical, and since these
facts are not essential in the proof of Theorem \ref{main}, we omit
them. The idea is to first to reduce to $k={\overline{\mathbb F}}_p$, and
then use the Galois action on  $k \otimes_{\mathbb F_p} \mathcal V$.

(v) For the same reason, we give no details here, and just point out
that the proof uses the fact that we have in reality the equality
${\rm rad} (k[H]) =  k \otimes_{\mathbb F_p} {\rm rad}\, ({\mathbb
F_ p}[H])$, which holds because this last ideal is clearly nilpotent.
The idea is to indicate that the invariants employed here are in
some sense defined over ${\mathbb F_p}$.
\end{proof}

\begin{defi}

  Let $\mathcal V$ a simple $\mathbb F_p [H]$-module.

  (i) Let $V \in \mathcal S _H$. We write $V|\mathcal V$ to say that $V$
  is a summand of $k \otimes_{\mathbb F_p} \mathcal V$.

 (ii) We will denote by $\deg \mathcal V$
the number of simple summands of $k \otimes_{\mathbb
  F_p} \mathcal V$ in a decomposition into simple $k[H]$-modules.
\end{defi}

\end{lem}

\subsubsection{Conclusion}

In view of Lemma \ref{rat} the two following lemmas end the proof of
Theorem \ref{main}.

\begin{lem}
Let $V\in \mathcal S _H$ such that $V|P$. Then :
$$\dim_{\mathbb F_p}{\rm Hom}_H(\pi_1^{ab}(Y),P)= \deg P\cdot
(h^2(H,V)+\delta_{Y,V})= \dim_{\mathbb F_p}
H^2(H,P)+\deg P\cdot\delta_{Y,V}$$ 
\end{lem}

\begin{proof}
  First notice that ${\rm Hom}(\pi_1^{ab}(Y),P) \simeq {\rm
  Hom}(\pi_1(Y),\mathbb F_p)\otimes_{\mathbb F_p} P$. But then using
  Artin-Schreier theory we see that :

  $${\rm Hom}(\pi_1(Y),\mathbb F_p) \simeq H^1(Y,\mathbb F_p) \simeq
  H^1(Y,\mathbb G_{a,Y})^F \simeq (H^0(Y, \Omega_Y)^{\mathcal
    C})^{\vee} $$

  As these isomorphisms commute with the $H$-action we conclude that :

  $$ {\rm Hom}_H(\pi_1^{ab}(Y), P) \simeq {\rm Hom}_H(H^0(Y,
\Omega_Y)^{\mathcal C},P)$$
  So it suffices to compute :
 $$\dim_{\mathbb F_p} {\rm Hom}_H (H^0(Y, \Omega_Y)^{\mathcal
    C},P) = \dim_k {\rm Hom}_H (H^0(Y, \Omega_Y)^s,
  k \otimes_{\mathbb F_p} P)$$

Because $k \otimes_{\mathbb F_p} P$ is semi-simple, we need only
to consider the largest semi-simple quotient of $H^0(Y, \Omega_Y)^s$,
that is $H^0(Y, \Omega_Y)^s / {\rm rad} (H^0(Y, \Omega_Y)^s)$. But
using Definition \ref{HW}, the isomorphisms
$\Omega^2 k / {\rm rad}\,\Omega^2 k \simeq P(\Omega^2 k) /\, {\rm rad}\,
P(\Omega^2 k)$ 
and \newline $P(\Omega^2 k) \simeq \mathop{\oplus}_{V' \in \mathcal S_H} P(V')
^ {\oplus h^2(H,V')}$ (see \cite{Kani}, Proposition 5), we finally see
that :
$$H^0(Y, \Omega_Y)^s / {\rm rad} (H^0(Y, \Omega_Y)^s) \simeq \mathop{\oplus}_{V' \in \mathcal S_H} {V'}
^{\oplus h^2 (H,V')+\delta_{Y,V'}}$$
hence :
$$\dim_k {\rm Hom}_H (H^0(Y, \Omega_Y)^s,
  k \otimes_{\mathbb F_p} P) = \sum _{V''|P}
  h^2(H,V'')+\delta_{Y,V''} = \deg P (h^2(H,V)+\delta_{Y,V}) $$

\end{proof}

\begin{lem}
Let $V\in \mathcal S _H$.
If the Kani inequalities hold, $G$ is a trivial extension of $H$ by $P$, and
$V|P$, then $\delta_{Y,V} >0$ 
\end{lem}

\begin{proof}
  Because of the Kani inequalities $\delta_{Y,V} \geq
  h^1(G,q^*V)-h^1(H,V)$, where $q$ is the epimorphism $q : G \thra H$. We
  will show that $h^1(G,q^*V)-h^1(H,V)=1$. In fact, using proposition
  \ref{infres}, and the hypothesis that $G$ is the trivial extension, we
  get an exact sequence :
  $$0 \ra H^1(H,V) \ra H^1(G, q^*V) \ra H^1(P, V)^H \ra 0$$

  But now : $$H^1(P, V)^H \simeq {\rm Hom}_{\mathbb F_p [H]}(P, V)\simeq
  {\rm Hom}_{k [H]}(k\otimes_{\mathbb F_p}P, V)\simeq k$$
\end{proof}

\section{Profinite proof for the extension of a $p'$-group by a $p$-group}

\label{prof}

In this section, we give another proof, almost purely algebraic,
of the main Theorem of \cite{PS} (i.e. Theorem \ref{main}
in the case of $H$ being a $p'$-group).

\subsection{Finite quotients of a free pro-$p$-group with $p'$-action}

Let $\gamma$ be a nonnegative integer and $F_\gamma$ be a free
pro-$p$-group with $\gamma$ generators, endowed with a continuous
action of a finite $p'$-group $H$. 

\begin{defi}
For each $\mathcal V \in \mathcal S _{\mathbb F_p [H]}$ let $\delta_V$
be the integer :
$$\delta_\mathcal V := \dim_{\mathbb F_p} H^1(F_\gamma \rtimes H, \mathcal V)$$ 
\end{defi}

\begin{rem}

In this definition, $\mathcal V$ is seen as a discrete $F_\gamma
\rtimes H$-module via the canonical map $F_\gamma \rtimes H \thra
H$. Using the exact sequence $0 \ra F_\gamma \ra  F_\gamma \rtimes H
\ra H \ra 0$ we can also consider $\mathcal V$ as a trivial continuous
$F_\gamma$-module, and the inflation-restriction sequence shows that
$\delta_\mathcal V = \dim_{\mathbb F_p} H^1(F_\gamma, \mathcal V)^H=
\dim_{\mathbb F_p}{\rm Hom}_H(F_\gamma, \mathcal V)$. But any continuous morphism $F_\gamma \ra \mathcal V$ is
trivial on the (profinite) Frattini subgroup $\Phi(F_\gamma)$, hence
$\delta_\mathcal V = \dim_{\mathbb F_p} {\rm
Hom}_H(F_\gamma/\Phi(F_\gamma), \mathcal V)$, which shows that the
$\delta_\mathcal V$'s are actually finite.

\end{rem}

\ \\
\begin{prop}
\label{alg}

Let $P$ a finite $p$-group endowed with an action of $H$.

The following assertions are equivalent :
\begin{itemize}
\item[(i)]
There exists a continuous $H$-epimorphism $ F_\gamma \thra P$ 
\item[(ii)]
$\forall \mathcal V \in \mathcal S _{\mathbb F_p [H]} \;\;\;
\dim_{\mathbb F_p} H^1(P \rtimes H, \mathcal V) \leq \delta_\mathcal
V$.
\end{itemize}
\end{prop}

\begin{proof}
It is quite clear that (i) implies (ii) : a $H$-epimorphism
$F_\gamma \thra P$ provides a continuous group epimorphism
$F_\gamma\rtimes H \thra P\rtimes H$, and the beginning of the
corresponding inflation-restriction sequence for $\mathcal  V$ :

$0 \ra H^1(P\rtimes H,\mathcal  V) \ra H^1(F_\gamma \rtimes H,\mathcal
V)\ra \cdots $

\noindent allows to conclude that (ii) holds.

To see that (ii) implies (i) first notice that, as in the remark
above,

\noindent
$\dim_{\mathbb F_p} H^1(P \rtimes H, \mathcal V)=\dim_{\mathbb F_p} {\rm
Hom}_H(P/\Phi(P), \mathcal V)$. Moreover, since $H$ is a $p'$-group,
the ring $\mathbb F_p [H]$ is semi-simple, and for any $\mathbb F_p
[H]$-module of finite type,
$\dim_{\mathbb F_p} {\rm Hom}_H(M, \mathcal V)$ is exactly $\deg \mathcal V$ times the
multiplicity of $\mathcal V$ in $M$. From this and (ii) we deduce the existence
of a surjective $H$-map $F_\gamma/\Phi(F_\gamma) \thra P/\Phi(P)$, and
so the existence of a group epimorphism $F_\gamma\rtimes H \thra
P/\Phi(P)\rtimes H$. Consider now the embedding problem~:

$$\xymatrix{ &&& F_\gamma\rtimes H \ar@{->>}[d]\ar@{-->}[dl] \\1
\ar[r] & \Phi(P)  \ar[r] & P\rtimes H \ar[r] & P/\Phi(P)\rtimes H \ar[r]  & 1  }  $$ 

According to \cite{Ser2}, Chapitre I, \S 4.2, Proposition 24, Corollaire 2, et
\S 3.3, Proposition 14, we certainly have ${\rm cd}_p( F_\gamma\rtimes
H) \leq 1$. Hence (see \ref{embpb})
this embedding problem has a weak solution $F_\gamma\rtimes H \ra
P\rtimes H$. This arrow commutes with projection on $H$, so we get a
similar diagram without $H$'s, but for continuous $H$-morphisms.
Now a classical argument (the Frattini subgroup $\Phi(P)$ is precisely
the set of non-generators in $P$) shows that the map $F_\gamma \thra
P$ must be an epimorphism. 
\end{proof}

\begin{rem}
One can of course ask if the $\delta_\mathcal V$'s do characterize the
action of $H$ on $F_\gamma$. We have no answer. 
\end{rem}

\subsection{Geometric application}

Let now $Y \ra X$ be a Galois \'etale cover of curves with Galois
group $H$, a $p'$-group. By pushing the exact sequence $0 \ra \pi_1(Y)
\ra \pi_1(X) \ra H \ra 0$ via the map $\pi_1(Y) \thra \pi_1^{(p)}(Y)$ we
obtain a split exact sequence (because ${\rm cd}_p H=0$)
which shows that $H$ acts naturally on
$\pi_1^{(p)}(Y)$. Since this group is known to be pro-$p$-free of rank
$\gamma_Y $ (the classical Hasse-Witt invariant of $Y$), we can
apply the result of the previous paragraph, with $\gamma = \gamma_Y$.

It is easy to see that Proposition \ref{alg} is equivalent to
Theorem \ref{main}. The main point is that $F_\gamma / \Phi(F_\gamma)
\simeq H^0(Y,\Omega_Y)^{\mathcal C}$, which allows to see that if $V\in
\mathcal S _{k[H]}$, $\mathcal V\in \mathcal S _{\mathbb F_p[H]}$ are
such that $V| \mathcal V$, then the $\delta_\mathcal V$'s of the previous
paragraph coincide with our former $\delta_{Y,V}$'s.

\end{document}